\documentclass{article}
\usepackage{quad}

\def\thetitle{Trees meet octahedron comparison}
\def\theauthors{Nina Lebedeva and Anton Petrunin}

\hypersetup{pdftitle={\thetitle},
pdfauthor={\theauthors}}

\newcommand{\Addresses}{{\bigskip\footnotesize

\noindent Nina Lebedeva,
\par\nopagebreak
 \textsc{St. Petersburg State University, 7/9 Universitetskaya nab., St. Petersburg, 199034, Russia}
\par
\nopagebreak
 \textsc{St. Petersburg Department of V. A. Steklov Institute of Mathematics of the Russian Academy of Sciences, 27 Fontanka nab., St. Petersburg, 191023, Russia}
  \par\nopagebreak
  \textit{Email}: \texttt{lebed@pdmi.ras.ru}

\medskip

\noindent   Anton Petrunin, 
\par\nopagebreak
 \textsc{Math. Dept. PSU, University Park, PA 16802, USA.}
  \par\nopagebreak
  \textit{Email}: \texttt{petrunin@math.psu.edu}
  
}}

\begin{document}

\title{\thetitle}
\author{\theauthors}

\date{}
\maketitle
\begin{abstract}
We show that metric trees and their products meet the octahedron comparison, which is a certain six-point metric comparison similar to Alexandrov's CAT(0) comparison. 
\end{abstract}

\parbf{Introduction.}
Let us recall the notion of graph comparison introduced in \cite{lebedeva-petrunin-zolotov}.

Suppose $\Gamma$ is a graph with vertices $v_1,\dots,v_n$.
A metric space $X$ meets the \emph{$\Gamma$-comparison} if for any set of points in $X$ labeled by vertices of $\Gamma$ there is a model configuration $\tilde v_1,\dots,\tilde v_n$ in the Hilbert space $\HH$ such that 
if $v_j$ is adjacent to $v_j$, then
$|\tilde v_i-\tilde v_j|_{\HH}\le | v_i-v_j|_{X}$,
and otherwise
$|\tilde v_i-\tilde v_j|_{\HH}\ge | v_i-v_j|_{X}$;
here $|p-q|_M$ denotes the distance between points $p$ and $q$ in the metric space~$M$.

\begin{wrapfigure}[8]{r}{25 mm}
\vskip-3mm
\centering
\includegraphics{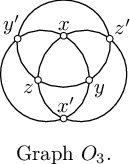}
\end{wrapfigure}

Let us denote by $O_3$ the graph of the octahedron.
For the labeling shown on the diagram, $O_3$-comparison means that any six-point configuration $x,x',y,y',z,z'\in X$ can be mapped to $\HH$ so that the diagonals $xx'$, $yy'$, and $zz'$ are not getting shorter
and the edges are not getting longer.

Note that the 4-cycle $C_4$ is an induced subgraph of $O_3$.
Therefore, \textit{$O_3$-comparison implies $C_4$-comparison};
the latter describes the $\CAT(0)$-comparison \cite{alexander-kapovitch-petrunin2022}.
In general, $C_4$-comparison does not imply $O_3$-comparison,
but for intrinsic metrics it is unknown.
Namely, the following question is open; it was motivated by \cite[1.19$_+(e)$]{gromov2007}.

\begin{thm}{Question}
Is it true that $O_3$-comparison holds in any geodesic $\CAT(0)$ space?
\end{thm}

An affirmative answer should lead to a classification of six-point spaces that admit an isometric embedding in a geodesic $\CAT(0)$ space.
On the other hand, as it follows from the next theorem,
a counterexample would mean that $O_3$-comparison is something new and interesting;
so, both answers are good.

Recall that a \emph{metric tree} is a geodesic space such that any two points $x$ and $y$ are connected by a unique geodesic $[xy]$,
and the union of any two geodesics $[xy]$, and $[yz]$ contains the geodesic $[xz]$.
Since metric trees and their products are $\CAT(0)$, the following theorem gives a partial answer to the question.

\begin{thm}{Theorem}
The $O_3$-comparison holds in products of metric trees.
\end{thm}

We suspect that our argument can be generalized to CAT(0) cube complexes and/or Euclidean buildings.
However, these generalizations are meaningful only if the question above has a negative answer.

Let us list the most relevant results.
\begin{itemize}
\item The so-called \emph{(4+2)-point comparison} is another six-point comparison;
it holds in any complete length $\CAT(0)$ space \cite{alexander-kapovitch-petrunin2011,alexander-kapovitch-petrunin2022}.
Some six-point metric spaces meet the (4+2)-point comparison, but do not admit a distance-preserving embedding into a complete length $\CAT(0)$
(see the discussion right after 7.2 in the updated arXiv version of \cite{alexander-kapovitch-petrunin2011}).
The $O_3$-comparison might be a sufficient condition for the existence of such embedding.

\item \textit{If a five-point space satisfies $C_4$-comparison, then it is isometric to a subset in a complete length $\CAT(0)$ space} \cite{toyoda2020,lebedeva-petrunin2021}.

\item \textit{For any integer $n\ge 4$, the $C_4$-comparison implies $C_n$-comparison} \cite{toyoda2023};
here $C_n$ denotes the $n$-cycle.
Here the metric does not have to be intrinsic. 
\end{itemize}
For other related questions, see \cite{gromov2001,gromov2007,petrunin2017,lebedeva,lebedeva-petrunin2010,lebedeva-petrunin2021,lebedeva-petrunin2022,lebedeva-petrunin2023,lebedeva-petrunin-zolotov, foertsch-lytchak-schroeder,eskenazis-mendel-naor,noar,toyoda2020,toyoda2023} and the references therein.

\bigskip

\parbf{Proof.}
First, let us show that \textit{$O_3$-comparison holds in any metric tree.}

Consider a metric tree $T$ with a six-point configuration $x,y,z,x',y',z'$ labeled by the vertices of $O_3$ on the diagram.
We can assume that the union of geodesics $[xx']$, $[yy']$, and $[zz']$ is connected.
Indeed, suppose $[xx']$ does not intersect $[yy']$ and $[zz']$.
Then we can shrink to a point a geodesic that minimizes the distance from $[xx']$ to $[yy']\cup[zz']$;
in the obtained tree $T'$ the distances $|x-x'|$, $|y-y'|$ and $|z-z'|$ do not change, and the remaining distances between $x,y,z,x',y',z'$ cannot get longer.
Therefore, any model configuration for the six points in $T'$ works for the original configuration in $T$.

Let 
\[Y=([xx']\cap [yy'])\cup([yy']\cap [zz'])\cup([zz']\cap [xx']).\]
Note that we are in one of the following two cases:
\begin{enumerate}
\item The set $Y$ is a \emph{tripod}; that is, $Y$ is a union of three geodesics meeting at one point, say $o$.
\item $Y$ lies in one of the  geodesics $[xx']$, $[yy']$, or $[zz']$.
\end{enumerate}

\begin{wrapfigure}{r}{35 mm}
\vskip-12mm
\centering
\includegraphics{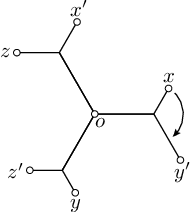}
\end{wrapfigure}

\parit{Case 1.}
Without loss of generality, we may assume that geodesics $[xx']$, $[yy']$, and $[zz']$ intersect as on the diagram (possibly with some degenerate edges).

Note that we can assume that in each pair $([ox],[oy'])$, $([oy],[oz'])$, and $([oz],[ox'])$ one of the geodesics lies in the other.
Indeed, assume it is not true, say $[ox]\not\subset[oy']$ and $|o-x|_T\le |o-y'|_T$.
Move $x$ to a point on $[oy']$ on the distance $|o-x|_T$ from $o$.
This way we decrease the distance $|x-y'|_T$, while the remaining fourteen distances between $x,y,z,x',y',z'$ do not change.
Therefore any model configuration for the new position of $x$ works for the old position.

If $[ox]\subset [oy']$, then define $a=x$, otherwise $a=y'$.
Similarly, choose  $b$ from $y$ or $z'$ and $c$ from $z$ or $x'$.
Let $[\tilde a\tilde b\tilde c]$ be a model plane for $[abc]$;
that is, $[\tilde a\tilde b\tilde c]$ is a plane triangle with the corresponding side lengths as in $[abc]$.

Note that $[ab]\subset [yy']$, $[bc]\subset [zz']$, and $[ac]\subset [xx']$.
Therefore we may choose points $\tilde x$, $\tilde x'$ on the line $\tilde a\tilde c$ so that the map $x\mapsto \tilde x$, $x'\mapsto \tilde x'$, $a\mapsto \tilde a$, $c\mapsto \tilde c$ is distance-preserving.
(The quadruple $\tilde x$, $\tilde x'$, $\tilde a$, $\tilde c$ might have identical points.)
Similarly, we may choose $\tilde y$, $\tilde y'$ on the line $\tilde a\tilde b$ and $\tilde z$, $\tilde z'$ on the line $\tilde b\tilde c$.
A couple of possible configurations are shown on the diagram.

\begin{figure}[ht!]
\centering
\includegraphics{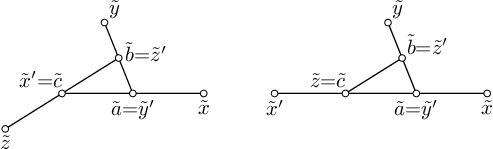}
\end{figure}

The triangle inequality implies that the obtained configuration meets the required conditions.

\parit{Case 2.}
Without loss of generality, we may assume that $Y\subset [xx']$.
Then $[yy']$ and $[zz']$ intersect $[xx']$ along two subsegments, possibly degenerate.
Without loss of generality, we may assume that geodesics $[xx']$, $[yy']$, and $[zz']$ are oriented the same way on their overlaps.

\begin{figure}[ht!]
\centering
\includegraphics{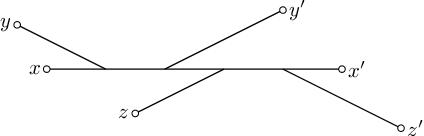}
\end{figure}

In this case, the needed configuration $\tilde x,\tilde x',\tilde y,\tilde y',\tilde z,\tilde z'$ on a line $\RR$ can be described by the coordinates:
\begin{align*}
\tilde x'&=|x'-x|_T,
&
\tilde y'&=|y'-x|_T,
&
\tilde z'&=|z'-x|_T,
\\
\tilde x&=0,
&
\tilde y&=|y'-x|_T-|y'-y|_T,
&
\tilde z&=|z'-x|_T-|z'-z|_T.
\end{align*}
The triangle inequality implies that this configuration meets the conditions. 

\parit{Final step.}
It remains to show that $O_3$-comparison holds for a product of trees $\Pi\z=T_1\times \dots\times T_n$.
Any six-point configuration $x,y,z,x',y',z'\in \Pi$ can be described by $n$ six-point configurations $x_i$, $y_i$, $z_i$, $x'_i$, $y'_i$, $z'_i\in T_i$ for each $i$.
Choose a model configuration $\tilde x_i$, $\tilde y_i$, $\tilde z_i$, $\tilde x'_i$, $\tilde y'_i$, $\tilde z'_i\in \HH$ for each $i$.
Observe that the $n$-tuples $\tilde x=(\tilde x_1,\dots, \tilde x_n),\z\dots,\tilde z'=(\tilde z_1,\dots, \tilde z_n)$ in $\HH=\HH^{\times n}$ form a model configuration of $x,y,z,x',y',z'$.
\qeds

\parbf{Acknowledgments.}
We want to thank Alexander Lytchak for his help. 

The first author was partially supported by the Russian Foundation for Basic Research grant 20-01-00070, the second author was partially supported by the National Science Foundation grant DMS-2005279
and the Ministry of Education and Science of the Russian Federation, grant 075-15-2022-289.

{\sloppy
\printbibliography[heading=bibintoc]
\fussy
}

\Addresses
\end{document}